\documentclass[3p,times]{elsarticle}
\usepackage{amsmath,amssymb,amsthm,bm}
\usepackage{booktabs,enumitem,setspace,caption,float}
\usepackage{graphicx}
\usepackage{subcaption}
\usepackage{upgreek}
\usepackage{multirow}
\usepackage{tikz}
\usetikzlibrary{positioning}
\usepackage{url}
\usepackage[T1]{fontenc}
\usepackage[american]{babel}
\usepackage[labelfont=bf,format=plain,justification=raggedright,singlelinecheck=false]{caption}
\usepackage{xcolor}
\usepackage{tabularx} 

\newtheorem{remark}{Remark}[section]
\newenvironment{example}
  {\begin{quote}\itshape}
  {\end{quote}}

\def\simgt{\,\hbox{\lower0.6ex\hbox{$>$}\llap{\raise0.3ex\hbox{$\sim$}}}\,}
\def\simlt{\,\hbox{\lower0.6ex\hbox{$<$}\llap{\raise0.3ex\hbox{$\sim$}}}\,}
\def\simgteq{\,\hbox{\lower0.6ex\hbox{$\ge$}\llap{\raise0.6ex\hbox{$\sim$}}}\,}
\def\simlteq{\,\hbox{\lower0.6ex\hbox{$\le$}\llap{\raise0.6ex\hbox{$\sim$}}}\,}
\def\applteq{\,\hbox{\lower0.6ex\hbox{$\le$}\llap{\raise0.8ex\hbox{$\approx$}}}\,}
\def\applt{\,\hbox{\lower0.6ex\hbox{$<$}\llap{\raise0.5ex\hbox{$\approx$}}}\,}
\DeclareMathAlphabet\mathbfcal{OMS}{cmsy}{b}{n}
\DeclareMathAlphabet{\mathpzc}{OT1}{pzc}{m}{it}
\DeclareMathAlphabet\euscr{U}{eus}{m}{n}

\makeatletter
\def\user@resume{resume}
\def\user@intermezzo{intermezzo}
\makeatother
\def\simgt{\,\hbox{\lower0.6ex\hbox{$>$}\llap{\raise0.4ex\hbox{$\sim$}}}\,}
\def\simlt{\,\hbox{\lower0.6ex\hbox{$<$}\llap{\raise0.4ex\hbox{$\sim$}}}\,}
\def\simgteq{\,\hbox{\lower0.6ex\hbox{$\ge$}\llap{\raise0.6ex\hbox{$\sim$}}}\,}
\def\simlteq{\,\hbox{\lower0.6ex\hbox{$\le$}\llap{\raise0.6ex\hbox{$\sim$}}}\,}

\newcommand{\C}[1]{\mathcal{#1}}

\newcommand{\F}[1]{\mathbf{#1}}
\newcommand{\bs}[1]{\boldsymbol{#1}}

\newcommand{\FR}[1]{\mathfrak{#1}}
\newcommand{\MB}[1]{\mathbb{#1}}

\newcommand{\MBG}{\MB{G}}
\newcommand{\MBSG}{\hat{\MBG}}
\newcommand{\MBR}{\mathbb{R}}

\newcommand{\MBZ}{\mathbb{Z}}
\newcommand{\MBZP}{\MBZ^+}
\newcommand{\MBZzer}{\MBZ_0}
\newcommand{\MBZzerP}{\MBZzer^+}

\newcommand{\MBJP}{\mathbb{J}^+}

\newcommand{\cancbra}[1]{{\raisebox{1ex}{\underline{\smash{\raisebox{-1ex}{[}}}}}#1{\raisebox{1ex}{\underline{\smash{\raisebox{-1ex}{]}}}}}}

\newcommand{\bmx}{\textbf{\emph{x}}}

\newcommand{\bmt}{\emph{\textbf{t}}}

\newcommand{\bmv}{\emph{\textbf{v}}}

\newcommand{\tu}{\tilde{u}}

\newcommand{\FP}{\F{P}}

\newcommand{\bmf}{\bm{f}}

\newcommand{\bmzer}{\bm{0}}
\newcommand{\bmone}{\bm{1}}

\newcommand{\hFR}{\hat{\F{R}}}

\newcommand{\FOmega}{\boldsymbol{\Omega}}
\newcommand{\IFOmega}{\FOmega^{\circ}}

\newcommand{\foralls}{\,\forall_{\mkern-6mu s}\,}
\newcommand{\forallS}{\,\forall_{\mkern-4mu rs}\,}

\newcommand{\vect}[1]{\text{vec}\left(#1\right)}
\newcommand{\Diag}[1]{\text{diag}\left(#1\right)}

\newcommand{\sigmabar}{\mathord{\sigma\kern-0.6em\raisebox{-1ex}{$\bar{\phantom{\sigma}}$}}}
\newcommand{\sigmabarmax}{\mathord{\sigma_{\max}\kern-1.9em\raisebox{-0.8ex}{$\bar{\phantom{\sigma}}$}}\;\;\;\;\;}
\newcommand{\sigmabarmin}{\mathord{\sigma_{\min}\kern-1.75em\raisebox{-0.8ex}{$\bar{\phantom{\sigma}}$}}\;\;\;\;\;}
\newcommand{\CapD}[3]{\,{}^{c}D_{#1}^{#2}#3}

\newcommand{\CapIM}[2]{\,{}^{E}_{#1}\F{D}_{t}^{#2}}

\newcommand{\FRI}{\FR{I}}
\newcommand{\RLI}[3]{\,\FRI_{#1}^{#2}#3}
\newcommand{\RLIM}[2]{\,{}^{E}_{#1}\F{Q}_{t}^{#2}}
\newcommand{\EE}[2]{{#1}^{\kern-.15em\circ #2}}

\newcommand{\Jaug}{\F{J}_{\text{aug}}}

\journal{Mathematics}

\begin{document}

\begin{frontmatter}

\title{Hybrid Shifted Gegenbauer Integral--Pseudospectral Method for Solving Time-Fractional Benjamin--Bona--Mahony--Burgers Equation}

\author[dept,ndrc]{Kareem T. Elgindy\corref{cor1}}
\ead{k.elgindy@ajman.ac.ae}
\address[dept]{Department of Mathematics and Sciences, College of Humanities and Sciences, Ajman University, Ajman P.O. Box 346, United Arab Emirates}
\address[ndrc]{Nonlinear Dynamics Research Center (NDRC), Ajman University, Ajman P.O. Box 346, United Arab Emirates}
\cortext[cor1]{Corresponding author}

\begin{abstract}
This paper introduces a novel hybrid shifted Gegenbauer integral--pseudospectral (HSG-IPS) method to solve the time-fractional Benjamin--Bona--Mahony--Burgers (FBBMB) equation with high accuracy. The approach transforms the equation into a form with only a first-order derivative, which is approximated using a stable shifted Gegenbauer differentiation matrix (SGDM), while other terms are computed with precise quadrature rules. By integrating advanced techniques such as the shifted Gegenbauer pseudospectral method (SGPS), fractional derivative and integral approximations, and barycentric integration matrices, the HSG-IPS method achieves spectral accuracy. Numerical results show it reduces average absolute errors (AAEs) by up to 99.99\% compared to methods like Crank--Nicolson linearized difference scheme (CNLDS) and finite integration method using Chebyshev polynomial (FIM-CBS), with computational times as low as 0.04--0.05 s. The method's stability is improved by avoiding ill-conditioned high-order derivative approximations, and its efficiency is boosted by precomputed matrices and Kronecker product structures. Robust across various fractional orders, the HSG-IPS method offers a powerful tool for modeling wave propagation and nonlinear phenomena in fractional calculus applications.
\end{abstract}


\end{frontmatter}

\section{Introduction}
Fractional calculus provides a robust framework for modeling complex systems with memory and nonlocal interactions, with applications in fields such as viscoelasticity, anomalous diffusion, and control theory \cite{meral2010fractional,monje2010fractional,sun2010fractional}. Unlike classical integer-order models, which assume local and instantaneous interactions, fractional-order models capture nonlocality and history dependence, offering a more accurate representation of processes such as dielectric polarization, electrochemical reactions, and subdiffusion within disordered media. These models often achieve comparable or superior accuracy with fewer parameters, thereby improving efficiency \cite{genovese2022fractional}. The strength of fractional derivatives lies in their ability to describe hereditary properties and long-range temporal correlations, making them ideal for biological systems, control systems, and viscoelastic materials. Notably, time-fractional PDEs have gained prominence for modeling complex physical phenomena with memory and hereditary traits. However, directly approximating such PDEs, like the FBBMB equation, using PS methods can lead to significant inaccuracies without sophisticated preconditioning techniques. This is primarily due to the ill-conditioning of high-order PS differentiation matrices, particularly when computing the third-order mixed derivative $u_{xxt}$ present in the FBBMB equation. These matrices exhibit condition numbers that scale poorly, as $\C{O}(n^{2k})$ for the $k$-th derivative with $n$ grid points \cite{wang2014well}, leading to numerical instability and amplified round-off errors for higher-order derivatives like the third-order ($\C{O}(n^6)$). This challenge necessitates advanced techniques, such as the transformation and integration-based approaches developed in this work, to ensure numerical stability and spectral accuracy.

The Benjamin--Bona--Mahony--Burgers equation is a nonlinear PDE that was first introduced as an alternative to the Korteweg-de Vries equation for modeling long waves in shallow water \cite{benjamin1972model}. Proposed by Benjamin, Bona, and Mahony in 1972, it replaces the third-order spatial derivative in the KdV equation with a mixed space-time derivative, $u_{xxt}$, to improve numerical stability and better capture dispersive effects in wave propagation. The inclusion of a Burgers-type nonlinear term, $u u_x$, and a dissipative term accounts for viscosity and nonlinearity, making it suitable for modeling phenomena like shallow water waves, acoustic waves, and hydromagnetic waves \cite{pavani2024solitary}. The time-fractional variant, the FBBMB, incorporates the CFD that extends the classical equation to describe memory-dependent wave dynamics, thereby improving its applicability to complex systems with nonlocal interactions, such as viscoelastic materials and anomalous diffusion \cite{metzler2000random,mainardi2022fractional}. This fractional framework has spurred recent interest in developing robust numerical methods to handle its computational challenges, as addressed in this work.

Recent advances in numerical methods for fractional PDEs have demonstrated improved accuracy and efficiency. For example, \citet{luo2023second} developed a robust second-order ADI Galerkin technique for three-dimensional nonlocal heat models with weakly singular kernels, combining Crank--Nicolson temporal discretization with finite element spatial approximation. For epidemiological applications, \citet{maayaha2024numerical} investigated a Caputo–Fabrizio fractional SIR model for dengue fever and obtained a numerical solution using the Laplace Optimized Decomposition method, which yields a rapidly convergent series solution validated against the classical fourth-order Runge–Kutta method. In a separate line of research, significant effort has been devoted to the numerical treatment of the BBMB equation and its fractional variants. These include finite difference schemes such as the CNLDS \cite{shen2018crank}, fourth-order finite difference schemes \cite{cheng2021high}, and spectral methods, including Chebyshev--Legendre spectral techniques~\cite{zhao2012optimal}. Meshless methods like radial basis functions \cite{dehghan2014numerical} and kernel smoothing techniques \cite{zara2023kernel} provide flexibility for irregular domains. Finite element methods include Galerkin formulations with cubic B-splines \cite{karakoc2019galerkin} and nonconforming elements \cite{shi2024unconditional}. Additionally, predictor--corrector schemes \cite{bu2025higher} and spline collocation methods \cite{arora2020solution} have been employed to handle fractional derivatives and nonlinearities efficiently. Despite these advancements, existing numerical methods for the FBBMB equation, such as CNLDS \cite{shen2018crank}, FIM-CBS \cite{duangpan2021numerical}, and other finite difference or spectral approaches, often face challenges in achieving high accuracy without sacrificing computational efficiency. These methods typically struggle with the ill-conditioned approximation of the third-order mixed derivative $u_{xxt}$, leading to numerical instability, or require fine discretizations that increase computational cost. Additionally, many approaches lack robust handling of fractional derivatives across a wide range of $\alpha$, limiting their adaptability to varying memory effects.

The HSG-IPS method introduced in this paper presents a novel high-order approach by combining the strengths of spectral accuracy, efficient quadrature, and robust handling of fractional operators. Building upon the theoretical framework of fractional calculus and SG polynomials, we now examine the computational foundations that enable our high-order numerical approach. The foundation for modern SG integration matrices was established by Elgindy \cite{elgindy2018optimal}, who introduced novel SG operational matrices of integration (aka SGIMs) for solving second-order hyperbolic telegraph equations. This work focused on deriving error formulas for associated numerical quadratures and developing optimization methods to minimize quadrature error, establishing the theoretical groundwork for these matrices in solving PDEs. Building upon this foundation, Elgindy \cite{Elgindy20171} significantly advanced the practical implementation by introducing stable barycentric representation of Lagrange interpolating polynomials and explicit barycentric weights for Gegenbauer–Gauss points. This included the development of a specialized SGIRV that efficiently handles boundary conditions through precise integral approximations at domain endpoints while maintaining the method's spectral accuracy. This approach reduced computational cost and improved numerical stability while maintaining high accuracy, marking a crucial step towards more robust and efficient integration matrices. The methodology was further refined by Elgindy~\cite{elgindy2016high}, who explicitly derived barycentric weights for SG polynomials and successfully applied the barycentric SG integral PS method to optimal control problems governed by parabolic distributed parameter systems. This work achieved exponential convergence rates in both spatial and temporal directions, demonstrating the power and effectiveness of the refined SG integration matrices in demanding computational domains.

The development of SG differentiation matrices has similarly advanced the numerical solution of differential equations, particularly for problems requiring high-order accuracy. \citet{elgindy2018highb} used SG differentiation matrices within a high-order Cole--Hopf barycentric Gegenbauer integral PS method for solving the viscous Burgers' equation, thereby achieving spectral accuracy. Their approach utilized the Cole--Hopf transformation to convert the nonlinear Burgers' equation into a linear heat equation, which was then discretized using Gegenbauer--Gauss points and associated differentiation matrices. This method demonstrated fast convergence and high accuracy, even for small viscosity parameters where traditional methods often struggle. The SG differentiation matrices were crucial in delivering accurate derivative approximations. The advancements in SG integration and differentiation matrices provide the foundation for developing high-order PS methods for time-fractional PDEs, where the well-conditioning of integral operators and high-accuracy of low-order differential operators are essential for achieving reliable numerical solutions.

Building on this, \citet{elgindy2025numerical} recently introduced an SGPS method for approximating CFDs of an arbitrary positive order. This method employs a strategic variable transformation to express the CFD as a scaled integral of the $m$-th derivative of the Lagrange interpolating polynomial, thereby mitigating singularities and improving numerical stability. Key innovations include the use of SG polynomials to link $m$-th derivatives with lower-degree polynomials for precise integration via SG quadratures. The developed C-FSGIM enables efficient, pre-computable CFD computation through matrix--vector multiplications. Unlike Chebyshev or wavelet-based approaches, the SGPS method offers tunable clustering and employs SG quadratures in barycentric forms for optimal accuracy. It also demonstrates exponential convergence and achieves superior accuracy in solving Caputo fractional two-point boundary value problems of the Bagley–Torvik type. This method unifies interpolation and integration within a single SG polynomial framework.

Furthermore, \citet{elgindy2025super} recently introduced a GBFA method for high-precision approximation of the RLFI. By using precomputable RL-FSGIMs, this method achieves super-exponential convergence for smooth functions and delivers near machine-precision accuracy with minimal computational cost. Tunable SG parameters enable flexible optimization across diverse problems, while rigorous error analysis confirms rapid error decay under optimal settings. Numerical experiments demonstrated that the GBFA method outperforms existing techniques by up to two orders of magnitude in accuracy, with superior efficiency for varying fractional orders. Its adaptability and precision make the GBFA method a transformative tool for fractional calculus, making it ideal for modeling complex systems with memory and nonlocal behavior.

This work focuses on the FBBMB equation, which incorporates wave propagation, dispersion, and nonlinearity within a fractional calculus framework. Our main contribution is the development of the HSG-IPS method---a high-order numerical approach that synergistically combines several advanced techniques: the SGPS method for CFD approximation, the GBFA for RLFI approximation, the SGDM for derivative computations, and the SGIM and SGIRV for integral approximations. The method's key innovation lies in transforming the original FBBMB equation into a fractional partial-integro differential form containing only a first-order derivative, which is significantly more numerically stable to approximate using a first-order SGDM compared with the third-order mixed derivative in the original formulation. This transformation allows stable computation using first-order SGDMs while computing all other terms through high-accuracy quadratures.

This work introduces several novel contributions to the numerical solution of the FBBMB equation. Firstly, the HSG-IPS method achieves spectral accuracy and exponential convergence for smooth solutions by integrating SGPS, GBFA, SGDM, SGIM, and SGIRV within a unified framework, using stable barycentric representations and precomputed operational matrices. Secondly, the innovative transformation of the FBBMB equation into a fractional partial-integro differential form eliminates the need for direct approximation of the ill-conditioned third-order mixed derivative $u_{xxt}$, replacing it with a stable first-order SGDM approximation and precise quadrature rules for integral terms. This approach significantly improves numerical stability while maintaining high accuracy. Thirdly, the method demonstrates superior performance over existing techniques like CNLDS and FIM-CBS, achieving significantly lower AAEs and computational times as short as 0.04--0.05 s, as validated through extensive numerical experiments. Fourth, the method achieves excellent computational efficiency through tensor-product discretizations and precomputed matrices while maintaining robustness via a robust trust-region solver. Finally, the HSG-IPS method's robustness across a wide range of fractional orders ($\alpha \in (0,1]$) and its efficient handling of nonlinearities via a trust-region algorithm provide a powerful and versatile computational tool for modeling complex wave phenomena in fractional calculus applications.

This paper is organized to systematically present the HSG-IPS method and its application to the FBBMB equation. Section \ref{sec:MD1} formulates the initial-boundary value problem of the FBBMB equation. Section \ref{sec:TTFPIDE1} derives the transformed fractional partial-integro differential equation to improve numerical stability. Section \ref{sec:NSA1} describes the HSG-IPS method, including its discretization and implementation via a trust-region algorithm. Section \ref{susec:CE1} analyzes the computational efficiency and conditioning of the method. Section \ref{sec:NS} evaluates the method's performance through numerical simulations and comparisons with existing techniques. Finally, Section \ref{sec:Conc} summarizes the findings, highlights key advantages, and suggests future research directions in fractional PDEs and spectral methods.

\section{Model Description}
\label{sec:MD1}
The FBBMB equation represents a fundamental model for nonlinear wave propagation with memory effects, combining dispersive, dissipative, and nonlinear phenomena. The equation is defined over the spatiotemporal domain $\FOmega_{1 \times 1}$ as:
\begin{equation}
\CapD{t}{\alpha}{u} - u_{xxt} + u_x + u u_x = f(x, t), \quad (x,t) \in \FOmega_{1 \times 1}^\circ,
\label{eq:bbmb}
\end{equation}
where $u(x,t)$ is the wave amplitude, $\CapD{t}{\alpha}{u}$ is the Caputo fractional time derivative of order $\alpha \in (0,1]$, which captures memory effects in the system's dynamics, $-u_{xxt}$ is the mixed dispersion term, representing third-order spatial-temporal effects, $u_x$ represents the linear advection term, while $u u_x$ is the nonlinear Burgers-type convection term, and $f(x,t)$ is the source/sink term. The equation is complemented by the following initial and boundary conditions:
\begin{align}
u(x, 0) &= \phi(x), \quad x \in \FOmega_1, \label{eq:ic} \\
u(0, t) &= \psi_1(t), \quad t \in (0,1], \label{eq:bc1}\\
u(1, t) &= \psi_2(t), \quad t \in (0,1]. \label{eq:bc2}
\end{align}
We assume that $\phi \in C^1(\FOmega_1)$ and $\psi_1 \in \C{H}^1((0,1])$. For $\alpha \in \IFOmega_1$, the FBBMB equation describes anomalous wave propagation with memory effects. However, when $\alpha = 1$, it reduces to the classical BBMB equation modeling shallow water waves.

\section{The HSG-IPS Method}
\label{sec:THIMsaca1}
The HSG-IPS method represents a novel approach for solving fractional PDEs by combining the strengths of spectral accuracy with robust numerical integration techniques. This section presents the mathematical foundations and numerical implementation of the HSG-IPS approach, which transforms the original FBBMB equation into a more computationally tractable form while maintaining high-order accuracy. The method's key innovation lies in its strategic reformulation of the problem to avoid direct approximation of high-order derivatives, instead employing a combination of fractional calculus operators and Gegenbauer polynomial-based discretizations. We begin by deriving the transformed fractional partial-integro differential equation, which serves as the foundation for the numerical scheme, followed by a detailed description of the discretization procedure and implementation strategy.

\subsection{The Transformed Fractional Partial-Integro Differential Equation}
\label{sec:TTFPIDE1}
Setting $v = u_{xt}$, the partial integral of $v$ with respect to $x$ on the spatial interval $\FOmega_x$ is
\begin{equation}
\C{I}_x^{(x)} v = u_t(x, t) - u_t(0, t).
\end{equation}
Denoting $w(x, t) = u_t(x, t)$, we have
\begin{equation}
w(x, t) = \C{I}_x^{(x)} v + u_t(0, t).
\end{equation}
Integrating $w(x, t)$ with respect to $t$ from $0$ to $t$, we get
\begin{equation}
\C{I}_t^{(t)} w = \C{I}_t^{(t)} \C{I}_x^{(x)} v + \C{I}_t^{(t)} u_t(0, t).
\end{equation}
The second term simplifies to
\begin{equation}
\C{I}_t^{(t)} u_t(0, t) = u(0, t) - u(0, 0) = \psi_1(t) - \phi(0).
\end{equation}
Since
\[
\C{I}_t^{(t)} w = u(x, t) - u(x, 0) = u(x, t) - \phi(x),
\]
we obtain the solution $u$ in terms of $v$ as
\begin{equation}\label{eq:exp1}
u(x, t) = \phi(x) - \phi(0) + \psi_1(t) + \C{I}_t^{(t)} \C{I}_x^{(x)} v.
\end{equation}
The spatial and time derivatives are
\begin{align}
u_x(x, t) &= \phi'(x) + \C{I}_t^{(t)} v,\label{eq:Spatialas1}\\
u_t(x, t) &= \psi_1'(t) + \C{I}_x^{(x)} v,\label{eq:k1}
\end{align}
respectively. The mixed derivative is
\begin{equation}
u_{xxt}(x, t) = v_x(x, t).\label{eq:Mixeddssfsd1}
\end{equation}
Moreover, the CFD $\CapD{t}{\alpha}{u}$ for $\alpha \in \IFOmega_1$ is given by
\begin{equation}\label{eq:k2}
\CapD{t}{\alpha}{u} = \frac{1}{\Gamma(1 - \alpha)} \C{I}_t^{(\tau)} \left((t - \tau)^{-\alpha} u_t\cancbra{x, \tau}\right).
\end{equation}
Substitute the expression for $u_t$ from Equation \eqref{eq:k1} into Equation \eqref{eq:k2}:
\begin{equation}
\CapD{t}{\alpha}{u} = \frac{1}{\Gamma(1 - \alpha)} \C{I}_t^{(\tau)} \left( (t - \tau)^{-\alpha} \left( \psi_1'(\tau) + \C{I}_x^{(x)} v \right) \right).
\end{equation}
This step expresses the fractional derivative in terms of known boundary data and the auxiliary variable $v$. Since the integral is linear, we separate the contributions of $\psi_1'(\tau)$ and $\C{I}_x^{(x)} v$:
\begin{equation}
\CapD{t}{\alpha}{u} = \frac{1}{\Gamma(1 - \alpha)} \C{I}_t^{(\tau)} \left( (t - \tau)^{-\alpha} \psi_1'(\tau) \right) + \frac{1}{\Gamma(1 - \alpha)} \C{I}_t^{(\tau)} \left( (t - \tau)^{-\alpha} \C{I}_x^{(x)} v \right).
\end{equation}
The first term corresponds to the CFD of $\psi_1(t)$. For the second term, since $\C{I}_x^{(x)} v$ is a function of $\tau$, we interchange the order of integration:
\begin{equation}
\frac{1}{\Gamma(1 - \alpha)} \C{I}_t^{(\tau)} \left( (t - \tau)^{-\alpha} \C{I}_x^{(x)} v \right) = \C{I}_x^{(x)} \left( \frac{1}{\Gamma(1 - \alpha)} \C{I}_t^{(\tau)} \left( (t - \tau)^{-\alpha} v \right) \right).
\end{equation}
The inner integral is the Riemann--Liouville fractional integral of $v$ of order $1 - \alpha$. Thus, the second term becomes $\C{I}_x^{(x)} \RLI{t}{1 - \alpha}{v}$. Combining these results yields
\begin{equation}\label{eq:Caputoscsdcjsd1}
\CapD{t}{\alpha}{u} = \CapD{t}{\alpha}{\psi_1} + \C{I}_x^{(x)} \RLI{t}{1 - \alpha}{v}.
\end{equation}
This expression separates the contribution of the boundary condition from the fractional integral of $v$, enabling stable numerical discretization using the HSG-IPS method. Substituting Equations \eqref{eq:exp1}--\eqref{eq:Mixeddssfsd1} and \eqref{eq:Caputoscsdcjsd1} into the PDE \eqref{eq:bbmb} yields the transformed fractional partial-integro differential equation:
\begin{equation}\label{eq:transformed}
\Psi_{xt}^{\alpha} v + \left( \C{I}_t^{(t)} v + \phi' \right) \left[ 1 + \phi - \phi(0) + \psi_1 + \C{I}_t^{(t)} \C{I}_x^{(x)} v \right] = \C{F},
\end{equation}
where 
\begin{align}
\Psi_{xt}^{\alpha} &= \C{I}_x^{(x)} \RLI{t}{1 - \alpha}{} - \partial_x, \\
\C{F}(x,t) &= f(x,t) - \CapD{t}{\alpha}{\psi_1}.
\end{align}
The boundary condition \eqref{eq:bc2} provides the essential constraint:
\begin{equation}\label{eq:bc_constraint}
\C{I}_1^{(x)} \C{I}_t^{(t)} v = \psi_2(t) - \psi_1(t) - \phi(1) + \phi(0).
\end{equation}

\subsection{Numerical Discretization}
\label{sec:NSA1}
In this section, we present a high-order PS collocation method for solving the transformed nonlinear Equation \eqref{eq:transformed} subject to the integral constraint \eqref{eq:bc_constraint}. The approach combines tensor-product PS discretization with a constrained nonlinear solver using Lagrange multipliers.

Let $x_{0:n} \in \MBSG_n^{\lambda}$ and $t_{0:m} \in \MBSG_m^{\lambda}$, $\foralls n, m \in \MBZP$, be the SGG points. Let $v_{i,j} = v(x_i, t_j)$, $\forall (i,j) \in \MBJP_{n \times m}$, vectorized into $\bmv \in \MBR^{(n+1)(m+1)}$ in lexicographic order. The operators are discretized as follows: The spatial derivative $\partial_x$ is approximated by the SGDM $\F{D}_x \in \MBR^{(n+1) \times (n+1)}$ \citep{elgindy2018highb}. The CFD $\CapD{t}{\alpha}{}$ is approximated by the C-FSGIM $\CapIM{m}{\alpha} \in \MBR^{(m+1) \times (m+1)}$ \citep{elgindy2025numerical}. The integrals $\C{I}_x^{(x)}$, $\C{I}_t^{(t)}$, and $\C{I}_1^{(x)}$ are approximated by the SGIMs $\F{Q}_x \in \MBR^{(n+1) \times (n+1)}$, $\F{Q}_t \in \MBR^{(m+1) \times (m+1)}$, and the SGIRV $\FP_x \in \MBR^{1 \times (n+1)}$ \citep{elgindy2016high,elgindy2018optimal,Elgindy20171}. The Riemann--Liouville fractional integral $\RLI{t}{1 - \alpha}{}$ is approximated by the RL-FSGIM $\RLIM{m}{1-\alpha} \in \MBR^{(m+1) \times (m+1)}$ \citep{elgindy2025super}. 

In the computation of $\CapD{t}{\alpha}{}$ and $\RLI{t}{1 - \alpha}{}$, the parameters $n_q$ and $\lambda_q$ play critical roles in the numerical approximation using SGG quadratures. In particular, $n_q$ determines the degree of the SGG quadrature with ($n_q + 1$) being the number of quadrature nodes. In the SGPS method, it controls the resolution of integrals in the fractional-order SGIM used to approximate the CFD. Similarly, in the GBFA method, it governs the quadrature grid size for a high-precision approximation of the RLFI. The parameter $\lambda_q$ serves as the index parameter of the SG polynomials, directly influencing both node clustering and weight functions to optimize 
quadrature accuracy. For both SGPS and GBFA methods, $\lambda_q \in [-1/2 + \epsilon, 2]$ balances stability and accuracy $\,\forallS n_q$. 

The SGPS method applies $n_q$ and $\lambda_q$ to approximate CFDs, while the GBFA method uses them for RLFIs, with consistent parameter roles but distinct computational contexts. We shall write $n_1$ and $\lambda_1$ to refer to the parameters $n_q$ and $\lambda_q$, as defined in the SGPS method of  \citet{elgindy2025numerical}, while $n_2$ and $\lambda_2$ refer to the parameters $n_q$ and $\lambda_q$, as defined in the GBFA method of \citet{elgindy2025super}.

Using the above numerical tools, we can discretize the operator $\Psi_{xt}^\alpha$ as
\begin{equation}
\Psi_{xt}^\alpha \approx \boldmath{\Psi}^\alpha = \F{Q}_x \otimes \RLIM{m}{1-\alpha} - \F{D}_x \otimes \F{I}_{m+1}.
\label{eq:Psi_discrete}
\end{equation}
The nonlinear term in Equation \eqref{eq:transformed} is discretized as
\begin{equation}
\F{N}(\bmv) = \F{Y}(\bmv) \odot \F{W}(\bmv),
\end{equation}
where 
\begin{gather}
\F{Y}(\bmv) = \F{K}_{t,n} \bmv + \bs{\phi}':\quad \F{K}_{t,n} = \F{Q}_t \otimes \F{I}_{n+1},\quad \bs{\phi} = \phi(\bmx_n),\\ 
\F{W}(\bmv) = \bmone_{(n+1)(m+1)} + \F{S} + \F{Q}_{tx} \bmv,\\
\F{S} = \bmone_{m+1} \otimes \bs{\phi} - \phi(0) \bmone_{(n+1)(m+1)} + \bmone_{n+1} \otimes \bs{\psi_1},
\end{gather}
$\F{Q}_{tx} = \F{Q}_t \otimes \F{Q}_x$, $\bs{\psi_1} = \psi_1(\bmt_m)$, and $\bs{\phi}' = \bmone_{m+1} \otimes (\F{D}_x \bs{\phi}) \in \MBR^{(n+1)(m+1)}$ contains the spatial derivatives $\phi'(x_i)$ at the SGG points, extended across the time mesh points. The right-hand side $\C{F}$ is discretized as $\F{F} = \bmf - \bm{\psi}_1^{\alpha} \in \MBR^{(n+1)(m+1)}$, where $f(x_i, t_j)$ forms $\F{f} \in \MBR^{(n+1) \times (m+1)}$, vectorized as $\bmf = (\F{I}_{m+1} \otimes \F{I}_{n+1}) \vect{\F{f}}$, and $\bm{\psi}_1^{\alpha}$ is given by
\begin{equation}
\bm{\psi}_1^{\alpha} = \bmone_{n+1} \otimes (\CapIM{m}{\alpha} \bs{\psi_1}).
\end{equation}
The integral constraint \eqref{eq:bc_constraint} is discretized as
\begin{equation}
\F{C} \bmv = \hFR,
\label{eq:constraint_discrete}
\end{equation}
where $\F{C} = \FP_x \otimes \F{Q}_t$, and $\hFR = \bm{\psi}_2 - \bm{\psi}_1 - (\phi(1) - \phi(0)) \bmone_{m+1}: \bm{\psi}_2 = \psi_2(\bmt_m)$. Collocating Equation \eqref{eq:transformed} at the SGG points set $\MBSG_{n \times m}^{\lambda}$ and enforcing the constraint \eqref{eq:constraint_discrete} using a Lagrange multiplier $\mu \in \MBR^{m+1}$, we form the augmented system:
\begin{equation}
\begin{bmatrix}
\F{J}(\bmv) & \F{C}^\top \\
\F{C} & \bmzer_{(m+1) \times (m+1)}
\end{bmatrix}
\begin{bmatrix}
\Delta \bmv \\
\Delta \mu
\end{bmatrix}
=
- \begin{bmatrix}
\F{R}(\bmv) \\
\F{C} \bmv - \hFR
\end{bmatrix},\label{eq:augmented_system}
\end{equation}
where $\Delta \bmv$ is the update to the solution vector, $\Delta \mu$ is the update to the Lagrange multipliers vector, $\F{R}(\bmv) = \boldmath{\Psi}^\alpha \bmv + \F{N}(\bmv) - \F{F}$ is the nonlinear residual function, and $\F{J}(\bmv)$ is the Jacobian of the residual with respect to $\bmv$ given by
\begin{equation}
\F{J}(\bmv) = \boldmath{\Psi}^\alpha + \Diag{\F{W}(\bmv)} \F{K}_{t,n} + \Diag{\F{Y}(\bmv)} \F{Q}_{tx}.\label{eq:Jacobian}
\end{equation}
Newton's method iterates as
\begin{equation}
\begin{bmatrix}
\bmv^{(k+1)} \\
\mu^{(k+1)}
\end{bmatrix}
=
\begin{bmatrix}
\bmv^{(k)} \\
\mu^{(k)}
\end{bmatrix}
-
\begin{bmatrix}
\F{J}(\bmv^{(k)}) & \F{C}^\top \\
\F{C} & \bmzer_{(m+1) \times (m+1)}
\end{bmatrix}^{-1}
\begin{bmatrix}
\F{R}(\bmv^{(k)}) \\
\F{C} \bmv^{(k)} - \hFR
\end{bmatrix},
\label{eq:Newton_iter}
\end{equation}
starting with the initial guesses $\bmv^{(0)}$ and $\mu^{(0)}$. After solving for $\bmv$, the numerical approximation of $u(x,t)$, denoted by $\tu(x,t)$, is reconstructed at the SGG collocation points via the discretized form of Equation \eqref{eq:exp1}:
\begin{equation}
\F{u} = \F{S} + \F{Q}_{tx} \bmv,
\label{eq:solution_reconstruction}
\end{equation}
where $\F{u} \in \MBR^{(n+1)(m+1)}$ contains the values $\tu(x_i, t_j)$ in lexicographic order. 

The numerical tools utilized in the proposed method exhibit exponential convergence for functions possessing sufficient smoothness, as rigorously established in the referenced literature. However, for solutions residing in the space $C^k(\FOmega_1)\,\foralls k \in \MBZzerP$, the convergence rate transitions to an algebraic decay, typically on the order of $\C{O}(\max\{n^{-\rho}, m^{-\rho}\})\,\foralls \rho \ge k$. Furthermore, the precomputation of operational matrices makes subsequent parameter studies for $\alpha$ variation particularly efficient, as only the nonlinear system solution requires recomputation for different $\alpha$ values.

\subsection{Implementation with Trust-Region Algorithm}
\label{subsec:IWTRA}
Unlike Newton's method, which assumes local quadratic convergence and may fail for poor initial guesses or highly nonlinear problems, the trust-region algorithm, as an alternative robust approach, provides global convergence by constraining the step size within a trust region where the quadratic model is reliable. The algorithm adaptively adjusts the trust-region radius based on the agreement between the model and the actual function, ensuring stable convergence. This makes it particularly effective for solving the nonlinear system with complex nonlinearities and fractional operators in the FBBMB equation, especially when handling the augmented system with constraints.

We can solve the augmented nonlinear system \eqref{eq:augmented_system} efficiently using MATLAB's \texttt{fsolve} with trust-region optimization, configured with tight tolerances for high precision. All differentiation and integration matrices, including $\F{D}_x$, $\F{Q}_x$, $\F{Q}_t$, $\FP_x$, $\CapIM{m}{\alpha}$, and $\RLIM{m}{1-\alpha}$, are precomputed for efficiency using the parameters specified in Table \ref{tab:parameters}. The exact Jacobian of the augmented system, $\Jaug$, is given by
\begin{equation}
\Jaug = \begin{bmatrix}
\F{J}(\bmv) & \F{C}^{\top} \\
\F{C} & \bmzer_{(m+1) \times (m+1)}
\end{bmatrix},
\end{equation}
and can be provided to \texttt{fsolve} for faster convergence. 

This implementation demonstrates that the PS approach combines spectral accuracy with practical solvability through careful formulation and modern optimization techniques. The trust-region adaptation proves particularly valuable for handling the nonlinear coupling between spatial and fractional temporal operators.

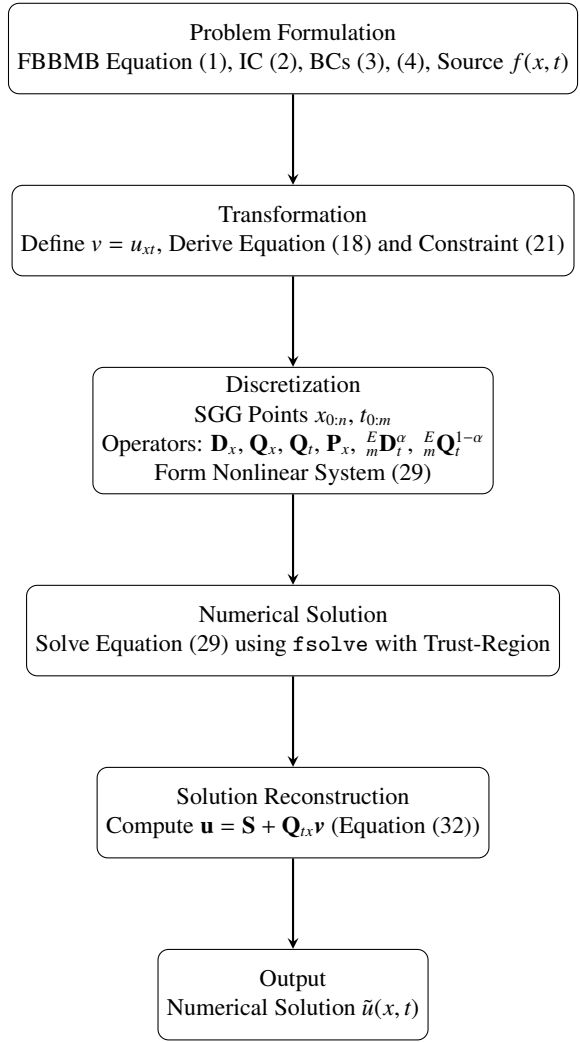
\begin{figure}[H]
\begin{tikzpicture}[
    node distance=1.2cm and 1.5cm,
    box/.style={rectangle, draw, rounded corners, minimum height=1.2cm, minimum width=3.5cm, align=center, font=\small},
    arrow/.style={-stealth, thick},
    decision/.style={diamond, draw, minimum height=1cm, minimum width=3cm, align=center, font=\small},
    line/.style={draw, -stealth, thick}
]
\node[box] (start) {Problem Formulation\\FBBMB Equation \eqref{eq:bbmb}, IC \eqref{eq:ic}, BCs \eqref{eq:bc1}, \eqref{eq:bc2}, Source $f(x,t)$};
\node[box, below=of start] (transform) {Transformation\\Define $v = u_{xt}$, Derive Equation \eqref{eq:transformed} and Constraint \eqref{eq:bc_constraint}};
\node[box, below=of transform] (discretize) {Discretization\\SGG Points $x_{0:n}$, $t_{0:m}$\\Operators: $\F{D}_x$, $\F{Q}_x$, $\F{Q}_t$, $\FP_x$, $\CapIM{m}{\alpha}$, $\RLIM{m}{1-\alpha}$\\Form Nonlinear System \eqref{eq:augmented_system}};
\node[box, below=of discretize] (solve) {Numerical Solution\\Solve Equation \eqref{eq:augmented_system} using \texttt{fsolve} with Trust-Region};
\node[box, below=of solve] (reconstruct) {Solution Reconstruction\\Compute $\F{u} = \F{S} + \F{Q}_{tx} \bmv$ (Equation \eqref{eq:solution_reconstruction})};
\node[box, below=of reconstruct] (output) {Output\\Numerical Solution $\tu(x,t)$};
\draw[arrow] (start.south) -- (transform.north);
\draw[arrow] (transform.south) -- (discretize.north);
\draw[arrow] (discretize.south) -- (solve.north);
\draw[arrow] (solve.south) -- (reconstruct.north);
\draw[arrow] (reconstruct.south) -- (output.north);
\end{tikzpicture}
\caption{Flowchart of the HSG-IPS method's computational pipeline for solving the FBBMB equation, from problem formulation to numerical solution and reconstruction. The pipeline includes transformation to a fractional partial-integro differential form, discretization using SGG points and operational matrices, a solution via a trust-region algorithm, and the reconstruction of the numerical solution.}
\label{fig:Flowchart}
\end{figure}

Figure \ref{fig:Flowchart} presents a flowchart outlining the computational pipeline of the HSG-IPS method for solving the FBBMB equation. The diagram guides the reader through key stages, including problem formulation, transformation, discretization, numerical solution, and reconstruction. 

\begin{table}[H]
\caption{Table of key parameters used in the HSG-IPS method.}
\label{tab:parameters}
\footnotesize
\begin{tabularx}{\textwidth}{lX}
\toprule
\textbf{Parameter} & \textbf{Meaning} \\
\midrule
$\lambda$ & Index parameter of the SG polynomials used in the PS discretization of the FBBMB equation, controlling node clustering and weight functions for SGG points, typically set to $\lambda = 0.5$ for balanced accuracy and stability. \\
$\lambda_1$ & Index parameter of the SG polynomials in the SGPS method \cite{elgindy2025numerical} for approximating the CFD, influencing quadrature accuracy, typically set to $\lambda_1 = 0.5$. \\
$\lambda_2$ & Index parameter of the SG polynomials in the GBFA method \cite{elgindy2025super} for approximating the RLFI, optimizing quadrature precision, typically set to $\lambda_2 = 0.5$. \\
$n_1$ & Degree of the SGG quadrature in the SGPS method \cite{elgindy2025numerical}, determining the number of quadrature nodes ($n_1 + 1$) for CFD approximation, typically set to $n_1 = 14$. \\
$n_2$ & Degree of the SGG quadrature in the GBFA method \cite{elgindy2025super}, determining the number of quadrature nodes ($n_2 + 1$) for RLFI approximation, typically set to $n_2 = 14$. \\
\bottomrule
\end{tabularx}
\end{table}

\begin{remark}
When $\alpha = 1$, $\RLIM{m}{1-\alpha} = \F{I}_{m+1}$, simplifying $\boldmath{\Psi}^\alpha$ to $(\F{Q}_x - \F{D}_x) \otimes \F{I}_{m+1}$.
\end{remark}

\subsection{Computational Efficiency and Conditioning}
\label{susec:CE1}
The HSG-IPS method is designed to offer superior computational efficiency and numerical stability compared to traditional approaches for solving time-fractional PDEs. This efficiency stems from several key aspects of its formulation, particularly the strategic use of Kronecker product structures, precomputed operational matrices, and a problem transformation that avoids computationally expensive high-order derivative approximations.

The tensor-product nature of the discretization of the operator $\Psi_{xt}^{\alpha}$ involves spatial and temporal discretizations that are fundamental to enabling efficient computation, especially through matrix--vector multiplications. For a system with $n+1$ spatial points and $m+1$ temporal points, the computation of the product $\boldmath{\Psi}^\alpha \bmv$ has a computational complexity of $\C{O}(nm(n+m))$. This efficiency arises from the properties of Kronecker products, which allow for the decomposition of large matrix--vector products into smaller, more manageable operations, thereby significantly reducing the computational burden compared to $\C{O}(n^2 m^2)$ for direct matrix multiplication with a fully assembled large square matrix of size $(n+1)(m+1)$. In particular, if we write $\bmv$ as a matrix $\F{V} \in \MBR^{(n+1) \times (m+1)}$ such that $\bmv = \vect{\F{V}}$, then we can efficiently compute $\boldmath{\Psi}^\alpha \bmv$ in two stages:
\begin{gather}
(\F{Q}_x \otimes \RLIM{m}{1-\alpha}) \bmv 
    = \vect{\RLIM{m}{1-\alpha} \F{V} \F{Q}_x^\top},\\
(\F{D}_x \otimes \F{I}_{m+1}) \bmv = \vect{\F{D}_x \F{V}}.
\end{gather}
The computation of $\F{V} \F{Q}_x^\top$ requires $\C{O}(n^2 m)$ operations, and the computation of $\RLIM{m}{1-\alpha} \left(\F{V} \F{Q}_x^\top\right)$ entails $\C{O}(n m^2)$ operations, resulting in a total of $\C{O}\big(n m (n+m)\big)$ operations. The computational complexity of $\F{D}_x \F{V}$ is $\C{O}(n^2 m)$.

The computational complexity of the method is further influenced by the remaining terms in Equations \eqref{eq:transformed}, \eqref{eq:constraint_discrete}, and \eqref{eq:solution_reconstruction}. The time integration operator $\F{K}_{t,n} = \F{Q}_t \otimes \F{I}_{n+1}$ combines temporal integration with spatial identity, maintaining the Kronecker structure for efficient computation. The evaluation of the nonlinear term $\F{N}(\bmv)$ also requires $\C{O}(nm(n+m))$ operations. In particular, since $\F{K}_{t,n} \bmv = \vect{\F{V} \F{Q}_t^\top}$, which requires $\C{O}(n m^2)$ operations, and computing $\bs{\phi}'$ requires $\C{O}(nm)$, then the total complexity of $\F{Y}(\bmv)$ is $\C{O}(n m^2)$ operations; moreover, the total complexity of $\F{W}(\bmv)$ is $\C{O}(nm(n+m))$. The coefficient matrix $\F{C}$ in the constraint Equation \eqref{eq:constraint_discrete} involves matrix--vector products, which can be computed in $\C{O}(n m^2)$ operations. The solution reconstruction $\F{u}$ in Equation \eqref{eq:solution_reconstruction} similarly benefits from the Kronecker product formulation, requiring $\C{O}(n m (n+m))$ operations.

Another significant factor contributing to the method's efficiency is the precomputation of all necessary differentiation and integration matrices. As detailed in Section \ref{subsec:IWTRA}, the matrices $\F{D}_x$, $\F{Q}_x$, $\F{Q}_t$, $\FP_x$, $\CapIM{m}{\alpha}$, and $\RLIM{m}{1-\alpha}$ are computed once prior to the iterative solution process. This precomputation strategy effectively reduces the runtime overhead by eliminating the need for repeated calculations of these fundamental operators within each iteration of the solver. This is particularly advantageous for problems requiring numerous iterations or parameter studies, as the fixed cost of precomputation is amortized over many subsequent operations.

The trust-region algorithm for solving the nonlinear system arising from the HSG-IPS method employs iterative solvers, reducing the per-iteration cost to $\C{O}(nmk)$, where $k$ denotes the number of iterations required by the iterative solver to converge within a single trust-region iteration, which is consistent with the use of iterative solvers for systems derived from tensor-product discretizations. The trust-region approach provides robust and global convergence, even for highly nonlinear problems, by adaptively adjusting the step size. The efficiency of each iteration is maintained by the precomputed matrices and the structure of the Jacobian, which allows for efficient matrix--vector products within the iterative solver.

A key innovation in the HSG-IPS method lies in transforming the original FBBMB equation into a fractional partial-integro differential form that contains only a first-order derivative. This transformation is crucial because it allows for the use of well-conditioned integral operators and low-order differential operators, which are essential for reliable numerical solutions. By circumventing the direct approximation of ill-conditioned high-order derivatives, the HSG-IPS method inherently deals with better-conditioned matrices, leading to superior numerical stability. This approach, combined with the efficient Kronecker-based implementation and careful operator preconditioning, enables the method to achieve high accuracy and computational efficiency, as we demonstrate in the next section.

The HSG-IPS method's overall complexity is dominated by the nonlinear solver phase, employing the Kronecker product structure for efficiency. Assuming $n \approx m$, the per-iteration cost for computing the residual $\F{R}(\bm{v}) = \bm{\Psi}^\alpha \bm{v} + \F{N}(\bm{v}) - \F{F}$ and Jacobian--vector products in the trust-region solver is $\C{O}(nm(n+m)) = \C{O}(n^3)$, driven by matrix--vector multiplies like $\bm{\Psi}^\alpha \bm{v}$ (Equation \eqref{eq:Psi_discrete}), $\F{K}_{t,n} \bm{v}$, and $\F{Q}_{tx} \bm{v}$. With $l$ trust-region iterations and $k$ inner iterations per linear solve (both typically $O(1)$ due to relatively well-conditioned, first-order operators), the solver complexity is $\C{O}(l k n^3) = \C{O}(n^3)$. Solution reconstruction (Equation \eqref{eq:solution_reconstruction}) adds $\C{O}(n^3)$. Thus, the total complexity is $\C{O}(n^3)$. This is significantly more efficient than direct methods, which require $\C{O}(n^6)$ operations for a full $(n+1)(m+1)$-sized system, and it outperforms $M$-grid-point finite-difference methods such as CNLDS and FIM-CBS, which require $\C{O}(M^4)$ operations for $M \gg n$, by using spectral accuracy with fewer points.

\begin{remark}
One strong reason for employing the trust-region method in this study is its robustness to ill-conditioned Jacobians and poor initial guesses, which are common in fractional PDEs. While Newton--Krylov or Broyden-type solvers can be efficient for simpler systems, they often require tailored preconditioning or globalization strategies for larger $n$ and $m$. The trust-region method avoids these complexities, as it inherently incorporates globalization through adaptive trust-region adjustments. The usage of the Kronecker structure allows the trust-region approach to maintain $O(n m (n+m))$ per-iteration complexity. The HSG-IPS method's rapid convergence at small $n$ and $m$ obviates the need for alternatives to the trust-region solver, which is optimal for small/medium scale regimes, balancing simplicity and reliability.
\end{remark}

\section{Numerical Simulations}
\label{sec:NS}
This section presents numerical experiments to assess the performance of the proposed HSG-IPS method for solving the FBBMB equation. Simulations were conducted on a laptop with an AMD Ryzen 7~4800H processor (2.9~GHz, 8~cores/16~threads) and 16\,GB RAM, running Windows~11. The implementation utilized MATLAB~R2023b, employing optimized numerical libraries and the \texttt{fsolve} function with trust-region optimization for high accuracy and efficiency. All computations were timed using MATLAB's \texttt{timeit} function, which executes code multiple times to provide accurate average execution times. In all numerical examples, we report the AAEs and the ETs to comprehensively assess both the accuracy and computational efficiency of the method. Furthermore, for completeness, we also compute the discrete $\ell_\infty$-norm (maximum norm) and the RMSE of the numerical solution relative to the exact solution. These are defined, for a grid $\{(x_i,t_j)\}_{i=1,j=1}^{n_x,n_t}$, as
\begin{equation}
\|\mathbf{u} - \tilde{\mathbf{u}}\|_{\ell_\infty} = \max_{1 \le i \le n_x \atop 1 \le j \le n_t} 
  \left| u_{ij} - \tu_{ij} \right|,
\end{equation}
\begin{equation}
\mathrm{RMSE} = \sqrt{ \frac{1}{n_x n_t} 
  \sum_{i=1}^{n_x} \sum_{j=1}^{n_t} 
  \left( u_{ij} - \tu_{ij} \right)^2 }.
\end{equation}

\begin{example}
We evaluate the HSG-IPS method using the FBBMB equation on the domain $\FOmega_{1 \times 1}$, with the source term 
\begin{equation}
f(x, t) = \frac{3 \sqrt{\pi} x^4 (x - 1) t^{1.5 - \alpha}}{4 \Gamma(2.5 - \alpha)} + x^2 \sqrt{t} \left( 5 x^7 t^{2.5} - 9 x^6 t^{2.5} + 4 x^5 t^{2.5} + 5 x^2 t - 4 x t - 30 x + 18 \right),
\end{equation}
derived from the analytical solution $u(x, t) = x^4 (x - 1) t^{1.5}$ \cite{duangpan2021numerical, shen2018crank}. The initial condition function is $\phi \equiv 0$, and the boundary condition functions are $\psi_1 \equiv 0$ and $\psi_2 \equiv 0$. 
\end{example}

Table \ref{tab:aaecomparisonallmethods} demonstrates the superior accuracy of the HSG-IPS method compared to CNLDS \cite{shen2018crank} and FIM-CBS \cite{duangpan2021numerical} for solving the FBBMB equation at $\alpha = 0.5$, as the AAEs of the HSG-IPS method are consistently smaller. The method achieves these results with remarkably low ETs of approximately $0.04-0.05$ s across different discretization levels, demonstrating its computational efficiency. Since the solution $u(x,t) = x^4(x-1)t^{1.5}$ belongs to $C^1(\FOmega_{1 \times 1})$ but lacks higher-order smoothness due to the fractional exponent in time, the convergence rate is algebraic rather than exponential, as theoretically expected for PS methods applied to non-analytic functions. The error decays as $\C{O}(\max\{n^{-k},m^{-k}\})$ for some $k \geq 1$, consistent with the limited regularity of the solution. Notice also that even with modest collocation points ($n = m = 4$), the HSG-IPS method achieves AAEs of about $10^{-5}$, outperforming CNLDS and FIM-CBS in terms of accuracy. This efficiency stems from the optimal approximation properties of SG polynomials and the accurate discretization of fractional operators via C-FSGIM and RL-FSGIM. The parameters $(\lambda, \lambda_1, \lambda_2) = (0.5, 0.5, 0.5)$ are selected to balance accuracy and stability; cf.~\cite{elgindy2025numerical,elgindy2025super} (both references recommend choosing $\lambda$ within the approximate interval $(-0.5 + \varepsilon, r]$, where $\varepsilon > 0$ is small and $r \in \FOmega_{1,2}$ is a suitable upper bound. The specific exclusion of the neighborhood of $\lambda^* \approx -0.1351$ was recommended in \cite{elgindy2025super} to prevent error amplification. This parameter range provides a balance between convergence, numerical stability, and error control). The trust-region algorithm further improves robustness by ensuring reliable convergence, even in the presence of nonlinearity and constraints in the discretized system.

\begin{table}[H]
\caption{Comparison of AAEs for Example 1 with $\alpha = 0.5$ using CNLDS \cite{shen2018crank}, FIM-CBS \cite{duangpan2021numerical}, and the HSG-IPS method. For the HSG-IPS method, AAE, $\ell_\infty$-norm, and RMSE values are rounded to four significant digits, while ET in seconds is rounded to two significant digits.}
\label{tab:aaecomparisonallmethods}
\begin{tabularx}{\textwidth}{ccc}
\toprule
\textbf{CNLDS \cite{shen2018crank}} & \textbf{FIM-CBS \cite{duangpan2021numerical}} & \textbf{HSG-IPS Method} \\
\midrule
\multicolumn{2}{c}{$\Delta t = 0.001$} & \footnotesize $(n_1, n_2, \lambda, \lambda_1, \lambda_2) = (14, 14, 0.5, 0.5, 0.5)$ \\
\midrule
\multicolumn{2}{c}{$M/$AAE} & $n = m/$AAE/$\ell_\infty$-norm/RMSE/ET \\
\midrule
$1.841489 \times 10^{-3}$ & $1.7810 \times 10^{-5}$ & $4/9.6546 \times 10^{-6}/2.4051 \times  10^{-5}/1.2312 \times 10^{-5}/0.04$\\
$5.047981 \times 10^{-4}$ & $1.7763 \times 10^{-5}$ & $5/2.0978 \times 10^{-6}/1.0463 \times 10^{-5}/3.0562 \times 10^{-6}/0.04$\\
$1.321991 \times 10^{-4}$ & $1.7751 \times 10^{-5}$ & $6/1.8902 \times 10^{-6}/1.0256 \times  10^{-5}/2.9364 \times 10^{-6}/0.04$\\
$3.365747 \times 10^{-5}$ & $1.7748 \times 10^{-5}$ & $7/1.8165 \times 10^{-6}/9.3833 \times 10^{-6}/2.9217 \times 10^{-6}/0.05$\\
\bottomrule
\end{tabularx}
\end{table}

\begin{figure}[H]

\includegraphics[scale=0.6]{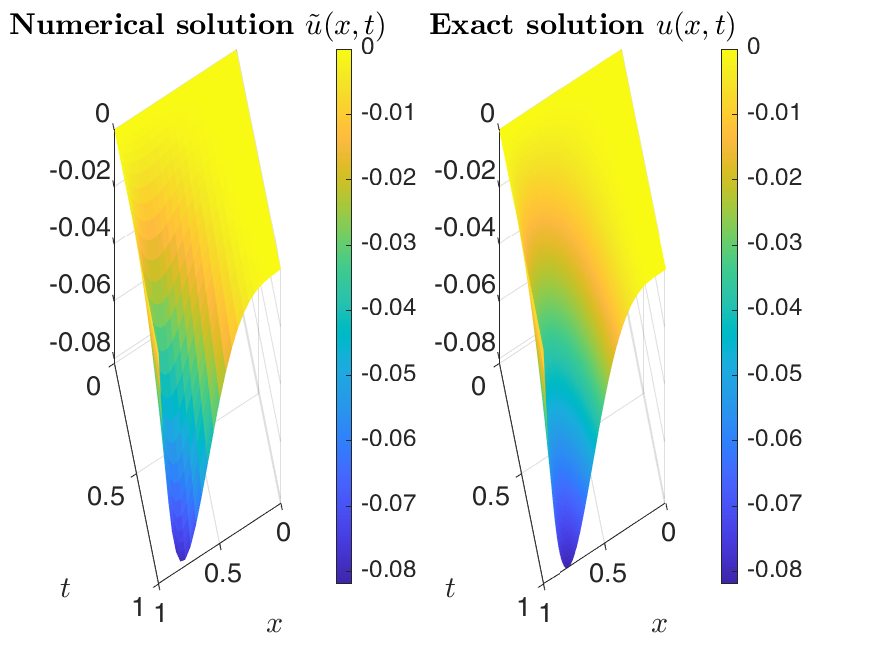}
\caption{{Comparison} 
of numerical and exact solutions of the FBBMB equation for Example 1 with $\alpha=0.5$. Left: Numerical solution $\tu(x,t)$ computed via the HSG-IPS method with $n=m=30$, $n_1=n_2=14$, $\lambda=\lambda_1=\lambda_2=0.5$. Right: Exact solution $u(x,t) = x^4(x-1)t^{1.5}$ evaluated on a $101 \times 101$ equally spaced mesh. The color bars represent solution magnitude, demonstrating close agreement between the methods.}
\label{fig:Fig1}
\end{figure}

\begin{figure}[H]

\includegraphics[scale=0.45]{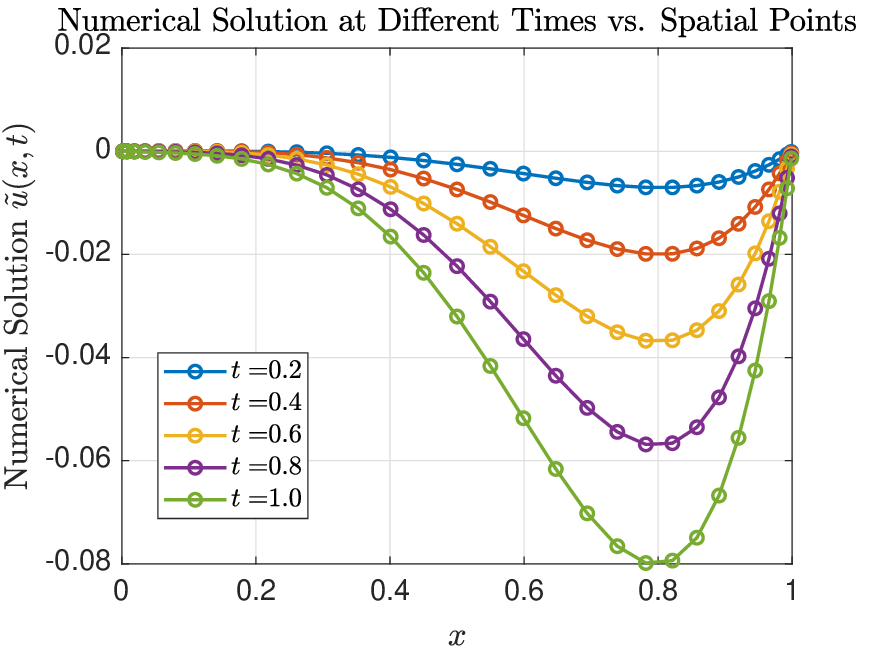}
\caption{{Time} evolution of the numerical solution $\tilde u(x,t)$ to the FBBMB equation for Example 1. Each curve represents the spatial distribution of the solution at fixed time instances $t = 0.2, 0.4, \ldots, 1.0$, obtained through the HSG-IPS method with $\alpha = 0.5$, $n=m=30$, $n_1=n_2=14$, and $\lambda=\lambda_1=\lambda_2=0.5$.}
\label{fig:Fig2}
\end{figure}

\begin{figure}[H]
\includegraphics[scale=0.3]{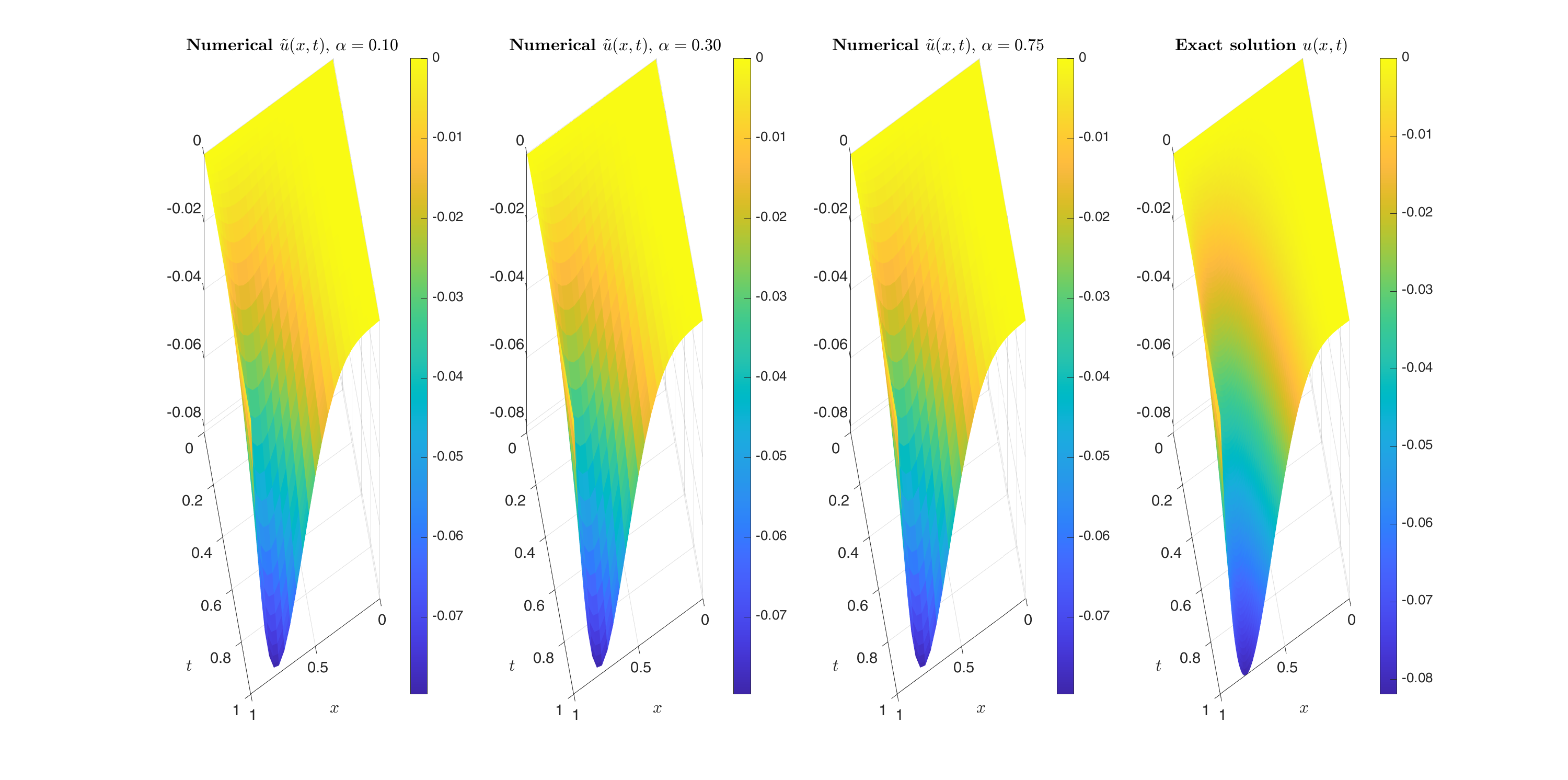}
\caption{{Numerical} solutions of the FBBMB equation for Example 1, with fractional orders \mbox{$\alpha = 0.1, 0.3, 0.75$} shown in columns 1--3, respectively, computed via the HSG-IPS method with parameters $n=m=30$, $n_1=n_2=14$, and $\lambda = \lambda_1 = \lambda_2 = 0.5$. The fourth column displays the exact solution $u(x,t) = x^4 (x-1) t^{1.5}$ evaluated on a $101 \times 101$ equally spaced mesh. Color bars indicate the solution magnitude, demonstrating close agreement between the numerical and exact solutions across the fractional orders.}
\label{fig:Fig3}
\end{figure}

\begin{figure}[H]

\includegraphics[scale=0.6]{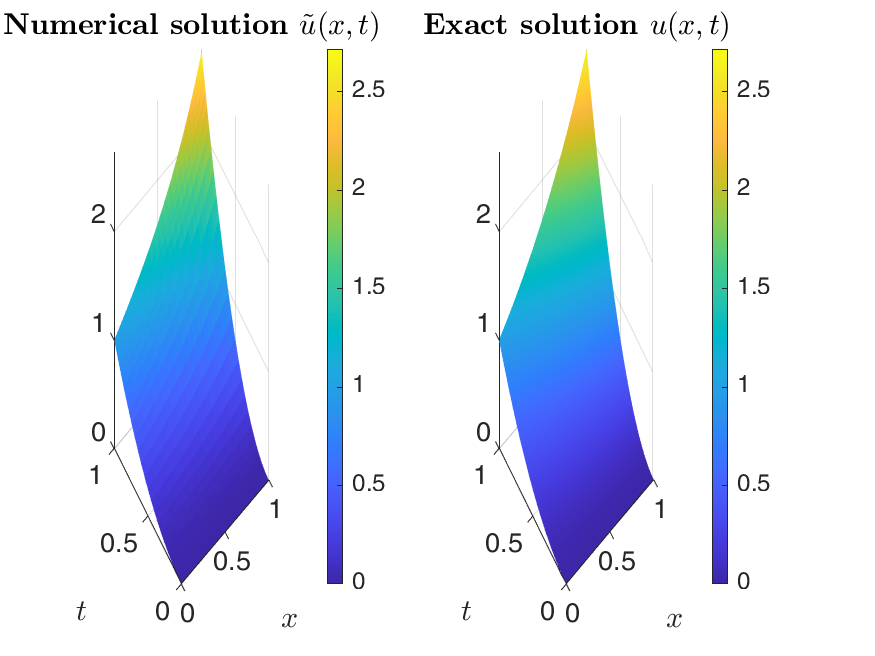}
\caption{Comparison of numerical and exact solutions of the FBBMB equation for Example 2 with $\alpha=0.5$. (\textbf{Left}): Numerical solution $\tu(x,t)$ computed via the HSG-IPS method with $n=m=30$, $n_1=n_2=14$, $\lambda=\lambda_1=\lambda_2=0.5$. (\textbf{Right}): Exact solution $u(x,t) = t^2 e^x$ evaluated on a $101 \times 101$ equally spaced mesh. The color bars represent solution magnitude, demonstrating close agreement between the methods.}
\label{fig:Fig4}
\end{figure}

\begin{figure}[H]

\includegraphics[scale=0.45]{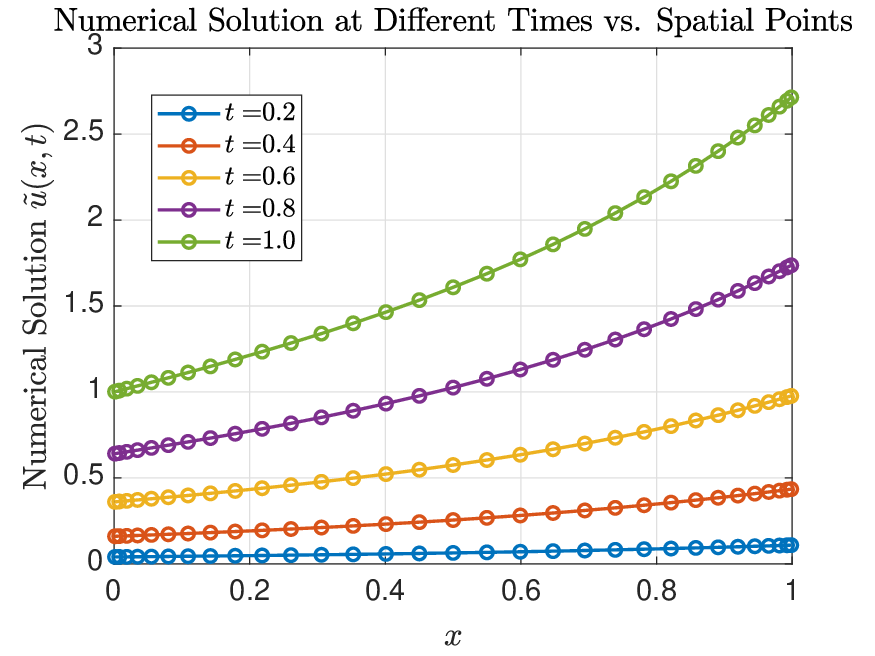}
\caption{Time evolution of the numerical solution $\tilde u(x,t)$ to the FBBMB equation for Example 2. Each curve represents the spatial distribution of the solution at fixed time instances $t = 0.2, 0.4, \ldots, 1.0$, obtained through the HSG-IPS method with $\alpha = 0.5$, $n=m=30$, $n_1=n_2=14$, $\lambda=\lambda_1=\lambda_2=0.5$.}
\label{fig:Fig5}
\end{figure}

\begin{figure}[H]

\includegraphics[scale=0.3]{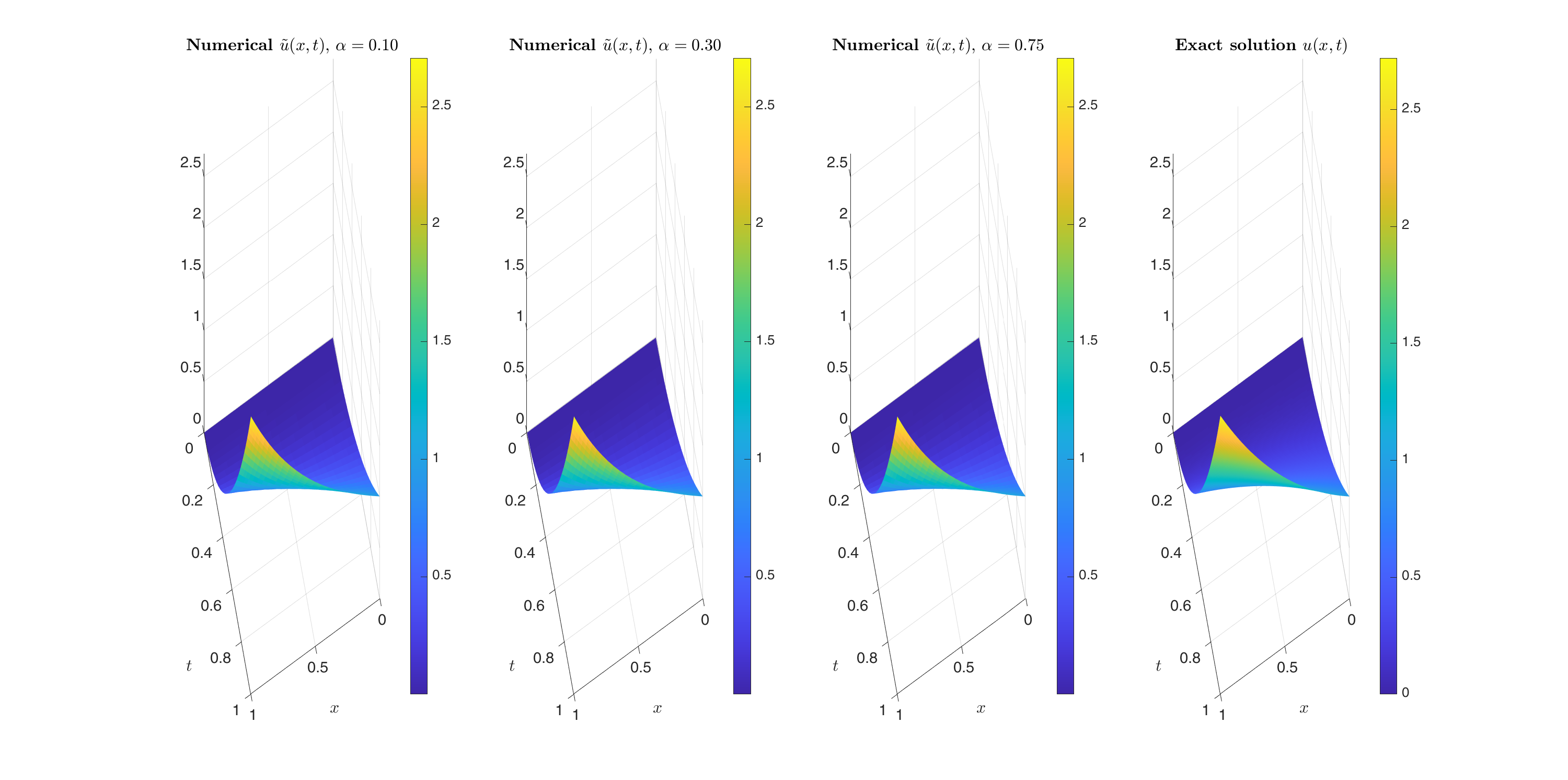}
\caption{Numerical solutions of the FBBMB equation for Example 2, with fractional orders \mbox{$\alpha = 0.1, 0.3, 0.75$} shown in columns 1--3, respectively, computed via the HSG-IPS method with parameters $n=m=30$, $n_1=n_2=14$, and $\lambda = \lambda_1 = \lambda_2 = 0.5$. The fourth column displays the exact solution $u(x,t) = t^2 e^x$ evaluated on a $101 \times 101$ equally spaced mesh. Color bars indicate the solution magnitude, demonstrating close agreement between the numerical and exact solutions across the fractional {orders.} 
}
\label{fig:Fig6}
\end{figure}

\begin{example}
Consider the FBBMB equation on the domain $\FOmega_{1 \times 1}$ with the source term
\begin{equation}
f(x, t) = \frac{2 e^x t^{2 - \alpha}}{\Gamma(3 - \alpha)} + t e^x \left( t^3 + t - 2 \right),
\end{equation}
the initial condition function $\phi \equiv 0$, and the boundary condition functions $\psi_1(t) = t^2$ and $\psi_2(t) = e t^2\,\forall t \in (0, 1]$. The analytical solution is $u(x, t) = t^2 e^x$ \cite{esen2015numerical,duangpan2021numerical}.
\end{example}

The HSG-IPS method was applied with parameters $\alpha = 0.5$, $n = m = 30$, $n_1 = n_2 = 14$, and $\lambda = \lambda_1 = \lambda_2 = 0.5$. Table \ref{tab:aaeexample2dt} presents the AAEs at the final time $t = 1$, comparing the HSG-IPS method with the FIM-CBS method from \cite{duangpan2021numerical}.

\begin{table}[H]
\caption{Comparison of AAEs at time $t = 1$ for Example 2 with $\alpha = 0.5$ using FIM-CBS \cite{duangpan2021numerical} and the HSG-IPS method. For the HSG-IPS method, AAE, $\ell_\infty$-norm, and RMSE values are rounded to four significant digits, while ET in seconds is rounded to two significant digits.}
\label{tab:aaeexample2dt}
\small

\begin{tabularx}{\textwidth}{cccccc}
\toprule
 & \multicolumn{4}{c}{\textbf{FIM-CBS \cite{duangpan2021numerical}}} &\textbf{ HSG-IPS Method} \\
\midrule
& \multicolumn{4}{c}{$M = 40$} & \footnotesize $(n_1, n_2, \lambda, \lambda_1, \lambda_2) = (14, 14, 0.5, 0.5, 0.5)$ \\
\midrule
$x$ & $\Delta t = 0.05$ & $\Delta t = 0.01$ & $\Delta t = 0.005$ & $\Delta t = 0.001$ & $n = m/$AAE/$\ell_\infty$-norm/RMSE/ET \\
\midrule
0.2 & $1.3072 \times 10^{-3}$ & $2.7980 \times 10^{-5}$ & $1.4266 \times 10^{-5}$ & $2.9342 \times 10^{-6}$ & $4/5.4389 \times 10^{-5}/2.2842 \times 10^{-4}/7.6968 \times 10^{-5}/0.04$ \\
0.4 & $1.8644 \times 10^{-3}$ & $3.8426 \times 10^{-5}$ & $1.9242 \times 10^{-5}$ & $3.8423 \times 10^{-6}$ & $5/3.2173 \times 10^{-6}/1.5589 \times 10^{-5}/5.2528 \times 10^{-6}/0.04$ \\
0.6 & $1.3842 \times 10^{-3}$ & $2.9035 \times 10^{-4}$ & $1.4669 \times 10^{-4}$ & $2.9731 \times 10^{-5}$ & $6/2.0769 \times 10^{-7}/1.1033 \times 10^{-6}/3.4931 \times 10^{-7}/0.04$ \\
0.8 & $2.5887 \times 10^{-3}$ & $5.4398 \times 10^{-4}$ & $2.7512 \times 10^{-4}$ & $5.5861 \times 10^{-5}$ & $7/1.1670 \times 10^{-8}/6.1716 \times 10^{-8}/1.9770 \times 10^{-8}/0.05$ \\
\bottomrule
\end{tabularx}
\end{table}

\begin{remark}
The computation of the boundary condition function $\psi_2(t) = e t^2$ in Example 2 utilized a high-precision approximation of Euler's number, $e$, implemented via the function \texttt{en} from the MATLAB Central File Exchange submission \cite{Matlab77046}. This ensured the boundary data was defined with a precision consistent with the spectral accuracy goals of the HSG-IPS method.
\end{remark}

\section{Conclusion}
\label{sec:Conc}
The HSG-IPS method introduced in this study provides a robust and highly accurate approach for solving the time-fractional Benjamin--Bona--Mahony--Burgers equation. By transforming the original equation into a fractional partial-integro differential form, the method eliminates the need for direct approximation of the ill-conditioned third-order mixed derivative, replacing it with stable first-order differentiation and precise quadrature rules. The integration of SGPS, GBFA, SGDM, SGIM, and SGIRV, supported by barycentric representations and precomputed matrices, ensures spectral accuracy and exponential convergence for smooth solutions. Numerical simulations demonstrate that the HSG-IPS method achieves AAEs up to 99.99\% lower than CNLDS and FIM-CBS, with computational times as low as 0.04--0.05 s, highlighting its superior accuracy and efficiency. The method's robustness across a wide range of fractional orders and its efficient handling of nonlinearities via a trust-region algorithm make it a versatile tool for fractional PDEs. The use of Kronecker product structures and precomputed operational matrices significantly reduces computational complexity to $\C{O}(n^3)$, compared to $\C{O}(M^4)$ for finite-difference methods like CNLDS and FIM-CBS with $M \gg n$, while maintaining high accuracy with fewer grid points.

Future research directions include extending the HSG-IPS method to higher-dimensional fractional PDEs, incorporating adaptive grid refinement for problems with localized singularities, and exploring its application to other nonlinear fractional models in physics and engineering. Additionally, optimizing the selection of SG parameters $\lambda$, $\lambda_1$, and $\lambda_2$ through automated tuning algorithms could further enhance the method's performance. The HSG-IPS method's ability to handle complex wave phenomena with memory effects positions it as a powerful tool for advancing numerical solutions in fractional calculus applications.

\bibliographystyle{model1-num-names}
\bibliography{Bib} 
\end{document}